\input ./style/arxiv-general.cfg
\documentclass[seceqn,MSNbibl,number,citesort,dvips]{arxbj}
\makeatletter
   \@ifpackageloaded{graphicx}{}{\usepackage{graphicx}}
\makeatother
\usepackage{upgreek}

\aid{0}
\volume{21}
\issue{3}
\pubyear{2015}
\firstpage{1824}
\lastpage{1843}
\doi{10.3150/14-BEJ627} 

\makeatletter

\newcommand{\vfrac}[2]{(#1)/#2}

\newcommand{\ifrac}[2]{(#1/#2)}

\newcommand{\rrvert}{\vert}
\newcommand{\llvert}{\vert}

\newtheorem{theorem}{Theorem}[section]

\newremark{rem}{Remark}[section]

\newtheorem{lemma}[theorem]{Lemma}

\newtheorem{prop}[theorem]{Proposition}

\newremark{Remark}{Remark}

\newcommand{\nn}{\nonumber}
\newcommand{\und}{\underline}
\newcommand{\La}{\Lambda}
\newcommand{\eps}{\varepsilon}
\newcommand{\ga}{\gamma}
\newcommand{\Om}{\Omega}
\newcommand{\si}{\sigma}
\newcommand{\eqref}[1]{(\ref{#1})}

\makeatother

\begin{document}
\begin{frontmatter}

\title{Extinction time for a random walk in a
random environment}
\runtitle{Extinction time}

\begin{aug}
\author[1]{\inits{A.}\fnms{Anna}~\snm{De Masi}\corref{}\thanksref{1}\ead[label=e1]{demasi@univaq.it}},
\author[2]{\inits{E.}\fnms{Errico}~\snm{Presutti}\thanksref{2}\ead[label=e2]{errico.presutti@gmail.com}},
\author[3]{\inits{D.}\fnms{Dimitrios}~\snm{Tsagkarogiannis}\thanksref{3}\ead[label=e3]{D.Tsagkarogiannis@sussex.ac.uk}} \and
\author[4]{\inits{M.E.}\fnms{Maria Eulalia}~\snm{Vares}\thanksref{4}\ead[label=e4]{eulalia@im.ufrj.br}}
\address[1]{Dipartimento di Ingegneria e Scienze dell'Informazione e
Matematica, Universit\`a di L'Aquila,
Via Vetoio, 1 67100 L'Aquila, Italy. \printead{e1}}
\address[2]{Gran Sasso Science Institute, Viale Francesco Crispi,7
67100 L'Aquila, Italy.\\ \printead{e2}}
\address[3]{Department of Mathematics,
University of Sussex,
Pevensey 2 Building,
Falmer Campus,
Brighton BN1 9QH, UK. \printead{e3}}
\address[4]{Instituto de Matem\'{a}tica, Universidade Federal do Rio de
Janeiro, Av. Athos da S. Ramos 149, 21941-909, Rio de Janeiro, RJ,
Brazil. \printead{e4}}
\end{aug}

\received{\smonth{3} \syear{2013}}
\revised{\smonth{12} \syear{2013}}

%
\begin{abstract}
We consider a random walk with death in $[-N,N]$ moving in a time dependent
environment. The environment is a system of particles which
describes
a current flux from $N$ to $-N$. Its evolution
is influenced by the presence of the random walk and in turn it affects
the jump rates of the random walk in a neighborhood of the endpoints,
determining also the rate for the random
walk to die. We prove an upper bound (uniform in $N$)
for the survival probability up to time $t$
which goes as $c\exp\{-b N^{-2} t\}$, with $c$ and $b$ positive constants.
\end{abstract}

%
\end{frontmatter}

\section{Introduction}
\label{sec:z1}

We consider a random walk on the discrete interval $\La_N:=[-N,N]$ of
$\mathbb Z$ which eventually dies by
jumping to a final state $\varnothing$ (where it stays thereafter).
Let $z\in\La_N \cup\{\varnothing\}$ denote the state of the random walk,
of which we say to be {\it alive} when $z\in\La_N$ and {\it dead}
when $z= \varnothing$.
When $z$ is alive and $|z|\le N-2$, it moves as a simple random walk:
after an
exponential time of mean 1 it jumps to its right or left neighbor with
probability $1/2$.
When $z\in I$, $I= I_+\cup I_-$, $I_+=\{N-1,N\}$, $I_-=\{-N,-N+1\}$ then,
besides moving, the walk $z$ may also die. The jump and death rates
depend on the environment.

The environment is a particle configuration $\eta$ on $\La_N
\setminus\{z\}$, $z$ the state
of the random walk (i.e., if $ z= \varnothing$ then $\eta\in\{0,1\}
^{\La_N}$, otherwise
$\eta\in\{0,1\}^{\La_N\setminus\{z\}}$). The evolution of the
environment is influenced by the motion of
the random walk: it consists of jumps of the particles (as second class
symmetric exclusion particles with $z$ being first class)
plus birth-death events localized in~$I$. The precise formulation is
given in the next section. We just
mention here that the birth-deaths events are ``rare''
as their intensity is proportional to $1/N$ and we are interested
in the case of large $N$.

When $z=\varnothing$, the environment evolves as
in \cite{DPTVjsp} with $K=2$ there ($K$ refers to the cardinality of
$I_+$ and $I_-$). Namely, it is the
simple symmetric exclusion process (SSEP, see \cite{li,li2})
in $\Lambda_N$ plus injection of particles into
$I_+$ and removal from $I_-$, from now on referred as the DPTV process:
at rate $j/(2N)$, one tries to inject a particle at the rightmost empty site
in $I_+$ and at the same rate there is an attempt to remove the
leftmost particle in $I_-$, the corresponding action being aborted if
$I_+$ is full or $I_-$ is empty. When $z\in\Lambda_N$, the evolution
of $(z,\eta)$ corresponds to a coupling of two realizations of the
DPTV process that differ at $z$ and undergo the same ``attempts'' to
create or to remove particles. Here, $j>0$ is a fixed parameter, while
we are interested in large $N$. Thus, when the random walk is dead,
the $\eta$ process describes a flux of particles from right to left
and it models
how currents can be induced by ``current reservoirs'', represented here
by the injection and removal processes at $I_+$ and, respectively, $I_-$.
With respect to \cite{DPTVjsp}, we now take $K=2$ for simplicity. It
will be clear that the arguments extend to any fixed $K$, the role
of current reservoirs being more closely achieved as $K$ grows.

The presence of the random walk
changes the picture and the purpose of this paper is to study how long
does such an influence
persist: we shall prove that the survival probability of the random walk
decreases exponentially in time, being bounded above by $c\exp\{-b
N^{-2}t\}$, $c,b>0$ independent of $t$ and $N$. In a companion paper
\cite{DPTVgap} we use the techniques and results developed here to
bound the extinction time in the case of several random walks. These
random walks correspond to the positions of discrepancies between two
configurations that evolve according to the DPTV process mentioned
before. By stochastic inequalities,
the result yields a lower bound of the form $b N^{-2}$ for the spectral
gap in this process, which is the motivation for our study here.


\section{Model and results}
\label{sec:z2}

The evolution of $(z,\eta)$ (random walk plus environment) is a Markov
process determined
by a generator $L$
which is the sum of the generators defined below, in \eqref
{z2.1}--\eqref{z2.8}. Letting the value
$\eta(x)=1$ ($\eta(x)=0)$ indicate the presence (absence) of a
particle at $x$, we may for convenience
always take $\eta\in\{0,1\}^{\Lambda_N}$ by requesting that $\eta
(z)=0$ whenever $z \ne\varnothing$.

We first suppose $z\ne\varnothing$ and write
%
\begin{eqnarray}
\label{z2.1} L^0_{\mathrm{env}} f(z,\eta)&=&\frac{1}2\Biggl\{
\sum_{x=-N}^{z-2} +\sum
_{x=z+1}^{N-1}\Biggr\}\bigl[f\bigl(z,\eta^{(x,x+1)}
\bigr)-f(z,\eta)\bigr],
\\
\label{z2.2} L^0_{\mathrm{z}} f(z,\eta)&=&\frac{1}2\bigl\{
\mathbf1_{z<N}\bigl[f\bigl(z+1,\eta ^{(z,z+1)}\bigr)-f(z,\eta)
\bigr]\nonumber\\[-8pt]\\[-8pt]
&&\hphantom{\frac{1}2\bigl\{}{}+\mathbf1_{z>-N}\bigl[f\bigl(z-1,\eta ^{(z-1,z)}\bigr)-f(z,
\eta)\bigr]\bigr\},\nonumber
\end{eqnarray}
where $\eta^{(x,x+1)}$ is obtained from $\eta$ by interchanging the
occupation values at $x$ and $x+1$,
and $\mathbf1_{z\in A}$ refers to the indicator function.

Denoting by
$\eta^{(+,x)}$ ($\eta^{(-,x)}$), the configuration which has the
value $1$ ($0$, resp.) at $x$ and otherwise coincides with $\eta$
%
\begin{eqnarray}
&&\hspace*{-21pt}\label{z2.3} L^+_{\rm env} f(z,\eta)  = \frac{j}{2N} \bigl\{
\mathbf1_{z<N} \bigl(1-\eta (N)\bigr) \bigl[f\bigl(z,\eta^{(+,N)}
\bigr) -f(z,\eta)\bigr]\nonumber
\\[-8pt]\\[-8pt]
&&\hspace*{-21pt}\hphantom{L^+_{\rm env} f(z,\eta)  =\frac{j}{2N} \bigl\{}{} +  \mathbf1_{z<N-1} \bigl(1-\eta(N-1)\bigr)\eta(N) \bigl[f\bigl(z,
\eta^{(+,N-1)}\bigr)-f(z,\eta )\bigr] \bigr\}\nn,
\\
&&\hspace*{-21pt}\label{z2.6} L^-_{\rm env} f(z,\eta)  = \frac{j}{2N} \bigl\{
\mathbf1_{z>-N} \eta (-N) \bigl[f\bigl(z,\eta^{(-,-N)}\bigr)-f(z,
\eta)\bigr]\nonumber
\\[-8.5pt]\\[-8.5pt]
 &&\hspace*{-21pt}\hphantom{L^-_{\rm env} f(z,\eta)  =\frac{j}{2N}
  \bigl\{}{}+  \mathbf1_{z>-N+1} \eta(-N+1) \bigl(1-\eta(-N)\bigr) \bigl[f\bigl(z,
\eta ^{(-,-N+1)}\bigr)-f(z,\eta)\bigr] \bigr\}\nn,
\\[-1pt]
&&\hspace*{-21pt}\label{z2.4} L^+_{\rm death} f(z,\eta)  = \frac{j}{2N} \bigl\{
\mathbf1_{z=N} \eta (N-1) \bigl[f\bigl(\varnothing,\eta^{(+,N)}
\bigr)-f(z,\eta)\bigr]\nonumber
\\[-8.5pt]\\[-8.5pt]
 &&\hspace*{-21pt}\hphantom{L^+_{\rm death} f(z,\eta)  =\frac{j}{2N} \bigl\{}{}+  \mathbf1_{z=N-1} \eta(N) \bigl[f\bigl(\varnothing,\eta^{(+,N-1)}
\bigr)-f(z,\eta)\bigr] \bigr\}\nn,
\\[-1pt]
&&\hspace*{-21pt}\label{z2.7} L^-_{\rm death} f(z,\eta)  = \frac{j}{2N} \bigl\{
\mathbf1_{z=-N} \bigl(1-\eta(-N+1)\bigr) \bigl[f(\varnothing,\eta)-f(z,\eta)\bigr]\nonumber
\\[-8.5pt]\\[-8.5pt]
&&\hspace*{-21pt} \hphantom{L^-_{\rm death} f(z,\eta)  =\frac{j}{2N} \bigl\{}{}+  \mathbf1_{z=-N+1} \bigl(1-\eta(-N)\bigr) \bigl[f(\varnothing,\eta)-f(z,
\eta)\bigr] \bigr\}\nn,
\\[-1pt]
&&\hspace*{-21pt}\label{z2.5} L^+_{\rm z} f(z,\eta)  = \frac{j}{2N}
\mathbf1_{z=N} \bigl(1-\eta(N-1)\bigr) \bigl[f\bigl(N-1,
\eta^{(+,N)}\bigr)-f(z,\eta)\bigr],
\\[-1pt]
&&\hspace*{-21pt}\label{z2.8} L^-_{\rm z} f(z,\eta)  = \frac{j}{2N}
\mathbf1_{z=-N} \eta(-N+1) \bigl[f\bigl(-N+1,\eta^{(-,-N+1)}
\bigr)-f(z,\eta)\bigr].
\end{eqnarray}
When $z=\varnothing$, the generator $L$ is the sum of only those in
\eqref{z2.1}, \eqref{z2.3} and \eqref{z2.6}
after replacing the indicator functions by 1 and putting $z=\varnothing$.
It is the one considered in \cite{DPTVjsp} in
the special case when the sets $I_{\pm}$ consist of only two sites.

Denote by $(z_t,\eta_t)_{t\ge0}$ the Markov process with the above
generator and by $P_{z,\eta}$ its law starting from $(z,\eta)$.
We now state the main result to be proven in the next sections.

\begin{theorem}
\label{thmz2.1}
There exist $c$ and $b$ positive and independent of $N$ so that for any
initial datum $(z_0,\eta_0)$, $z_0\ne\varnothing$ and any $t>0$
%
\begin{equation}\label{z2.9}
P_{z_0,\eta_0} [ z_t \ne\varnothing ] \le c \mathrm{e}^{-b N^{-2}t}.
\vspace*{-1pt}
\end{equation}
\end{theorem}

\section{The auxiliary process}
\label{sec:z3}

It will be useful to consider an auxiliary process $(\tilde z_t)_{t\ge
0}$. This will be a time-inhomogeneous Markov process whose
jump intensities at time $t$ are obtained by averaging those of the
original process over the environment conditioned on the state
of the random walk at that time. The explicit expression of the time
dependent generator $\mathcal L_t$ is given below
in \eqref{z3.3} after introducing some definitions and notation. We
fix hereafter arbitrarily
the initial condition $(z_0,\eta_0)$ at time $0$, $z_0\ne\varnothing$,
and denote by $\tilde P_{z_0}$ and $\tilde E_{z_0}$
the law of the auxiliary process and corresponding expectation. We
shall prove that for any bounded measurable
function $\phi(z,\eta)= f(z)$:
%
\begin{equation}\label{z3.0}
E_{z_0,\eta_0} \bigl[ \phi(z_t,\eta_t) \bigr] =
\tilde E_{z_0} \bigl[ f(\tilde z_t) \bigr] .
\end{equation}
By taking $f(z)=\mathbf1_{z\ne\varnothing}$, \eqref{z3.0} shows that the
distributions of the extinction time for the true and the auxiliary
processes are the same.
The proof of \eqref{z3.0} follows from the equality
%
\begin{equation}\label{z3.4}
\frac{\mathrm{d}}{\mathrm{d}t} E_{z_0,\eta_0} \bigl[ \phi(z_t,
\eta_t) \bigr] = E_{z_0,\eta_0} \bigl[\mathcal L_t
f(z_t) \bigr],
\end{equation}
which we shall prove next.\vadjust{\goodbreak}

We obviously have $L^{\pm}_{\rm env}\phi=0$ and, for $z\ne\varnothing$,
$L^0_{\rm z} \phi= \mathcal L^0 f$ with $\mathcal L^0$ the generator
of the
simple random walk on $[-N,N]$ with jumps outside $[-N,N]$ suppressed
(as in the definition of $L^0_{\rm z}$).
Recalling \eqref{z2.4}--\eqref{z2.7}
\begin{eqnarray*}
L^+_{\rm death} \phi&=& \frac{j}{2N} \bigl\{\mathbf1_{z=N}
\eta(N-1) \bigl[f(\varnothing)-f(N)\bigr] + \mathbf1_{z=N-1} \eta(N) \bigl[f(
\varnothing)-f(N-1)\bigr] \bigr\},
\\
L^-_{\rm death} \phi&=&\frac{j}{2N} \bigl\{\mathbf1_{z=-N}
\bigl(1-\eta (-N+1)\bigr) \bigl[f(\varnothing)-f(-N)\bigr]
\\
& &\hphantom{\frac{j}{2N} \bigl\{} {}+ \mathbf1_{z=-N+1} \bigl(1-\eta(-N)\bigr) \bigl[f(\varnothing)-f(-N+1)
\bigr] \bigr\}.
\end{eqnarray*}
By \eqref{z2.5} and \eqref{z2.8},
\begin{eqnarray*}
L^+_{\rm z} \phi&=&\frac{j}{2N} \mathbf1_{z=N} \bigl(1-
\eta(N-1)\bigr) \bigl[f(N-1)-f(N)\bigr],
\\
L^-_{\rm z} \phi&=&\frac{j}{2N} \mathbf1_{z=-N} \eta(-N+1)
\bigl[f(-N+1)-f(-N)\bigr].
\end{eqnarray*}
Thus, we define
%
\begin{eqnarray}\label{z3.1}
d(N,t) &=& \frac{j}{2N} E_{z_0,\eta_0}\bigl[\eta_t(N-1) |
z_t=N\bigr],\nn
\\
d(N-1,t) &=& \frac{j}{2N} E_{z_0,\eta_0}\bigl[\eta_t(N) |
z_t=N-1\bigr]\nn,
\\[-8pt]\\[-8pt]
d(-N,t) &=& \frac{j}{2N} E_{z_0,\eta_0}\bigl[\bigl(1-
\eta_t(-N+1)\bigr) | z_t=-N\bigr],\nn
\\
d(-N+1,t) &=& \frac{j}{2N} E_{z_0,\eta_0}\bigl[\bigl(1-
\eta_t(-N)\bigr) | z_t=-N+1\bigr] \nonumber
\end{eqnarray}
set $d(z,t)=0$ if $|z|<N-1$, and let
%
\begin{eqnarray}\label{z3.2}
a(N,t) &=& \frac{j}{2N} E_{z_0,\eta_0}\bigl[ \bigl(1-
\eta_t(N-1)\bigr) | z_t=N\bigr],\nn
\\[-8pt]\\[-8pt]
a(-N,t) &=& \frac{j}{2N} E_{z_0,\eta_0}\bigl[ \eta_t(-N+1) |
z_t=-N\bigr]. \nonumber
\end{eqnarray}
Given $t\ge0$, define
%
\begin{eqnarray}\label{z3.3'}
\mathcal L^a_tf(z) &=& \mathbf1_{z\neq\varnothing}\mathcal
L^0 f (z)+ \mathbf1_{z=N}a (N,t) \bigl[f(N-1)-f(N)\bigr]\nn
\\[-8pt]\\[-8pt]
&&{}+ \mathbf1_{z=-N}a (-N,t) \bigl[f(-N+1)-f(-N)\bigr] \nonumber
\end{eqnarray}
and
%
\begin{equation}\label{z3.3}
\mathcal L_tf(z)  =   \mathcal L^a_tf(z)+
d(z,t) \bigl[f(\varnothing)-f(z)\bigr],
\end{equation}
so that we get \eqref{z3.4}, and hence \eqref{z3.0} at once.

The auxiliary process $\tilde z_t$ is thus the Markov process with time
dependent generator
$\mathcal L_t$. It is a simple random walk with extra jumps $N\to N-1$
and $-N\to-N+1$
which occur with intensities $a(\pm N,t)$ and death rates ($z\to
\varnothing$) given by
$d(z,t)$.
Calling $\mathcal P_{z_0}$ the law of the process $\tilde z_t$ with
time dependent generator $\mathcal{L}^a_t$ (same fixed $\eta_0$ and
the same initial condition $z_0$ at time $0$) and denoting by $\mathcal
E_{z_0}$ the corresponding expectation, one sees that (see \cite{BG},
Chapter III),
%
\begin{eqnarray}\label{z3.5}
\nn P_{z_0,\eta_0} [ z_t \ne\varnothing ] &=& \tilde
P_{z_0} [ \tilde z_t \ne\varnothing ] = \mathcal
E_{z_0} \biggl[ \exp\biggl\{ - \int_0^t
d(\tilde z_s,s)\,\mathrm{d}s \biggr\} \biggr]
\\[-8pt]\\[-8pt]
&\le& \mathcal E_{z_0} \biggl[ \exp\biggl\{ - \int
_0^t d(N,s) \mathbf1_{\tilde z_s=N}\,\mathrm{d}s \biggr
\} \biggr] \nonumber,
\end{eqnarray}
where the last inequality is not really necessary, brings some loss,
but is just to simplify.

The proof of Theorem \ref{thmz2.1} follows from \eqref{z3.5} and the
following two statements which will be proved in the next sections.
\begin{itemize}
\item There are $\delta^*>0$ and $\kappa>0$ so that for all $t\ge
T_2=\kappa N^2$:
%
\begin{equation}\label{z3.6}
d(N,t)\ge\frac{j\delta^*}N.
\end{equation}

\item There exists a positive constant $b$ so that calling $T^*(t)$
the total time spent at $N$ by $\tilde z_s, 0 \le t$:
%
\begin{equation}\label{z3.7}
\mathcal E_{z_0} \bigl[ \mathrm{e}^{-j\delta^* N^{-1} T^*(t)} \bigr] \le \mathrm{e}^{-b
N^{-2} t},\quad\quad
t\ge T_2= \kappa N^2.
\end{equation}
\end{itemize}

\section{Proof of \texorpdfstring{\protect\eqref{z3.7}}{(3.9)}}
\label{sec:4}

Throughout the rest of the paper we shall write $\eps\equiv N^{-1}$.
With the notation introduced above and writing $\mathcal E_{t,\tilde
z}$ for the
conditional distribution (under $\mathcal P_{z_0}$) of $(\tilde z_s, s
\ge t)$ given $\tilde z_t=\tilde z$, we prove: Given any $\delta>0$
there is $p<1$ so that uniformly in $\eps$ and for all non negative
integers $n$:
%
\begin{equation}\label{z4.1}
\mathcal E_{t_n,\tilde z_{t_n}} \bigl[ \mathrm{e}^{-X}\bigr] \le p,\quad\quad
 X:= \eps\delta
\int_{t_n}^{t_{n+1}} \mathbf1_{\tilde z_s= N}\,\mathrm{d}s,\quad\quad
t_n= 2\eps^{-2} n.
\end{equation}
We see that \eqref{z3.7} follows at once from \eqref{z4.1}: taking
$\delta=j\delta^*$ in the latter and using the Markov property the
left-hand side of \eqref{z3.7} is bounded from above by $p^{[t/(2N^2)]}$,
compatible with its right-hand side. Now, the key point in proving \eqref{z4.1}
is the following.

\begin{lemma}
For any $0<c_-<c$, there is $p<1$ (as given in \eqref{z4.4} below) so
that the following holds.
Let $(\Om,\mu)$ be a probability space, $E$ the expectation and
$\mathcal F$ the set of all measurable functions $f\ge0$ such that
$E[f]\ge c_-$ and $E[f^2]\le c^2$. Then
$E[\mathrm{e}^{-f}] \le p$ for any $f\in\mathcal F$.
\end{lemma}

\begin{pf} Let $f\in\mathcal F$,
$\zeta:= c_-/2$, $\ga:= \mu[f> \zeta]$. Then
%
\begin{equation}\label{z4.2}
c_- \le E [ f] =E [ f; f \le\zeta]+ E [ f; f> \zeta] \le\zeta (1-\ga) + c
\ga^{1/2}.
\end{equation}
Call $a= \ga^{1/2}$, then \eqref{z4.2} yields
$\zeta(1-a^2) +ca -c_- \ge0$, so that $a_- <a<a_+$ where $a_{\pm}$
are the roots of the corresponding equation with equality:
\[
\zeta a^2 -ca+c_-- \zeta=0,\quad\quad \mbox{that is, } 2\zeta a= c \pm\sqrt{
c^2 - 4\zeta(c_--\zeta)} = c \pm\sqrt{c^2-c_-^2}.
\]
Thus,
\[
2\zeta a_- = c - c\sqrt{ 1 -\frac{c_-^2}{c^2}}\ge c - c \biggl( 1 -
\frac{1}2 \frac{c_-^2}{c^2} \biggr) = \frac{c_-^2}{2c}
\]
so that (since $\mu[f> \zeta] = a^2$ and $a\ge a_-$)
%
\begin{equation}\label{z4.3}
\mu[f> \zeta] \ge \biggl( \frac{c_-}{2c} \biggr)^2
\end{equation}
and
%
\begin{eqnarray}\label{z4.4}
E\bigl[\mathrm{e}^{-f}\bigr] &\le& \mathrm{e}^{-\zeta}\mu[ f> \zeta] +1- \mu[ f>
\zeta]\nn
\\[-8pt]\\[-8pt]
&= & 1-\mu[ f> \zeta]\bigl(1-\mathrm{e}^{-\zeta}\bigr) \le1- \biggl(
\frac{c_-}{2c} \biggr)^2 \bigl(1-\mathrm{e}^{-c_-/2}\bigr)=:p.\nonumber
\end{eqnarray}
\upqed\end{pf}

To apply the lemma, we need to prove the existence of constants
$0<c_-<c$ so that for any $\eps$, any $n$ and $\tilde z_{t_n}$,
%
\begin{equation}\label{z4.5}
c_-\le\mathcal E_{t_n,\tilde z_{t_n}} [ X ], \quad\quad\mathcal E_{t_n,\tilde
z_{t_{n}}} \bigl[
X^2 \bigr] \le c^2.
\end{equation}

Proof that $\mathcal E_{t_n,\tilde z_{t_n}} [ X ] \ge c_-$.
We claim that under $\mathcal P_{t_n,\tilde z_{t_n}}$ the time spent at
$N$ by the process $(\tilde z_t)$ during the time interval
$[t_n,t_{n+1}]$ is
stochastically larger than the time spent at $N$ during the interval
$[0, 2N^2]$ by a simple random walk $(x_t)$ in $\mathbb Z$ that starts
at time $0$ from $\tilde z_{t_n}$. Since $a(N,t) <1/2$, the intensity
with which the process $(\tilde z_t)$ jumps from $N$ to $N-1$ is
smaller than one, which is the jump rate of $(x_t)$. It is then easy to
construct a coupling of both processes for which $|x_{t-t_n}-N| \ge
|\tilde z_t -N|$ for all $t$. This is
done by constructing a suitable time inhomogeneous Markov process for
the pair evolution. Here are the details of the coupling, setting the
jump rates at time $t$ when $(x_{t-t_n}=x, \tilde z_t=\tilde z)$ with
the property that $|\tilde z-N|\le|x-N|$:
\begin{itemize}
\item
Let $x=\tilde z=N$, the pair $(x,\tilde z)$ moves to $(N-1,N-1)$ with
intensity $1/2$; it moves to $(N+1,N)$ with intensity $1/2-a(N,t)$, and
with intensity $a(N,t)$ it moves to $(N+1,N-1)$.

\item Let $\tilde z=N$ and $x\ne N$. From $(x,N)$, the pair moves
to $(x,N-1)$ with intensity $1/2+ a(N,t)$; with intensity $1/2$ it
moves to $(x-1,N)$ and with intensity $1/2$ it moves to $(x+1,N)$.

\item Let $\tilde z\ne\pm N$ and $x:|\tilde z-N| \le|x-N|$. With
intensity $1/2$ both coordinates move by 1 away from $N$, and with
intensity $1/2$ both move by 1 toward $N$.

\item Let $\tilde z= -N$ and $x:|\tilde z-N| \le|x-N|$. With
intensity $1/2$ $x$ moves by $1$ toward $N$ and $\tilde z$ moves to
$-N+1$; with intensity $a(-N,t)$, $x$ moves by $1$ away from $N$ and
$z$ moves to $-N+1$; with intensity $1/2-a(-N,t)$, $x$ moves by 1 away
from $N$ and $\tilde z$ stays put.
\end{itemize}

 Observe that this gives a coupling of the processes and that
all the jumps preserve the inequality $|\tilde z-N| \le|x-N|$.

Proof that $\mathcal E_{t_n, \tilde z_{t_n}} [ X^2 ] \le c^2$.
Since $\mathcal E_{t_n,\tilde z_{t_n}} [ X^2 ] \le\mathcal E_{t_n, N}
[ X^2 ]$, we just need to prove the inequality when $\tilde z_{t_n}=N$.
A simple construction, similar to the previous one, allows to couple
$(\tilde z_{t})$ and $(x_t)$ a simple random walk that moves in
$[0,N]$, that is, the jumps to $-1$ and $N+1$ are suppressed, starting
at $N$ at time $0$, in such a way that $\tilde z_{t} \le x_{t-t_n}$ for
all $t \in[t_n,t_{n+1}]$. The details are quite simple and, therefore,
omitted. As a consequence, the time spent at $N$ by $(\tilde z_t)$
during $[t_n,t_{n+1}]$ is stochastically smaller than that spent at $N$
during $[0, 2N^2]$ by this simple random walk $(x_t)$.

The process $(x_t)$ can be realized on the unit rate symmetric simple
random walk $(y_t)$ on $\mathbb Z$ (jumps $\pm1$ with rate $1/2$ each)
by identifying sites on $\mathbb Z$ modulo repeated reflections around
$N+1/2$ and $-1/2$, that is, reflections that identify $N+1$ with $N$, and
  $-1$ with $0$ (see, e.g., \cite{DPTVjsp}, Proposition 4.1). Thus,
calling $N_i$ the images of $N$ under the above reflections, we have to bound
%
\begin{equation}\label{z4.6}
2 \int_0^{t_1} \mathrm{d}s \int_s^{t_1} \mathrm{d}s' \sum_{i,k} E_{N} [
\mathbf1_{y_s=N_i} \mathbf1_{y_{s'}=N_k} ].
\end{equation}
By the local central limit theorem as in \cite{lawler} (see also
Theorem 3 in \cite{DPTVjsp}), this can be bounded in terms of Gaussian
integrals, from which
\eqref{z4.5} is proved. Details are omitted.

\section{Proof of \texorpdfstring{\protect\eqref{z3.6}}{(3.8)}}

We continue to write $\eps:= N^{-1}$, and set the following notation:
\[
\pi(x,t) := P_{z_0,\eta_0}[z_t=x]= \tilde P_{z_0}[
\tilde z_t=x],\quad\quad B(x,t):= (j\eps)^{-1}d(x,t)\pi(x,t),
\]
so that \eqref{z3.6} is implied by
%
\begin{equation}\label{z5.1}
B(N, t)\ge\delta^* \pi(N,t), \quad\quad t\ge T_2= \kappa\eps^{-2}.
\end{equation}
Having defined
%
\begin{equation}
\label{z5.2} T_1= \eps^{-(1-a)},\quad\quad T_0 =
T_1-\eps^{-(1-a)/2},\quad\quad T_2= \kappa\eps
^{-2}, \quad\quad\mbox{$a>0$ small enough}
\end{equation}
and
%
\begin{equation}\label{z5.3}
p_t(x,y) = \mbox{transition probability of the simple random walk on
$\Lambda_N$}
\end{equation}
(the jumps to $\pm(N+1)$ being suppressed),
we postpone the proof of the following three bounds, for $t \geq T_2$:
\begin{itemize}
\item There are $b_1>0$ and, for any $n$, $c_n$ so that
%
\begin{equation}
\label{z5.4} B(N,t)\ge b_1\sum_z
p_{T_1}(N,z)\pi(z,t-T_1)-c_n\eps^n
\tilde P_{z_0}[\tilde z_{t-T_2} \ne\varnothing].
\end{equation}

\item There are $b_2>0$, and for any $n$, $c_n$ so that
%
\begin{equation}
\label{z5.5} \pi(N,t)\le b_2 \sum_z
p_{T_1}(N,z)\pi(z,t-T_1)+c_n\eps^n
\tilde P_{z_0}[\tilde z_{t-T_2} \ne\varnothing].
\end{equation}
\item There is $b_3>0$ so that
%
\begin{equation}
\label{z5.6} \pi(N,t)\ge b_3 \eps^3 \tilde
P_{z_0}[\tilde z_{t-T_2} \ne\varnothing].
\end{equation}
\end{itemize}

\begin{Claim*}
\eqref{z5.1} follows from \eqref{z5.4}, \eqref{z5.5},
\eqref{z5.6}.
\end{Claim*}

\begin{pf}
By \eqref{z5.6}, we get from \eqref{z5.5}
%
\begin{equation}
\label{z5.7} \biggl[1-\frac{c_n}{b_3}\eps^{n-3}\biggr]\pi(N,t)\le
b_2 \sum_z p_{T_1}(N,z)\pi
(z,t-T_1),
\end{equation}
and from \eqref{z5.4}
%
\begin{equation}
\label{z5.8} B(N,t)\ge b_1\sum_z
p_{T_1}(N,z)\pi(z,t-T_1)-\frac{c_n}{b_3}\eps
^{n-3}\pi(N,t).
\end{equation}
Using \eqref{z5.7} and \eqref{z5.8}, we have
%
\begin{equation}
\label{z5.9} B(N,t)\ge\frac{b_1}{b_2} \biggl[1-\frac{c_n}{b_3}
\eps^{n-3}\biggr]\pi (N,t)-\frac{c_n}{b_3}\eps^{n-3}\pi(N,t),
\end{equation}
which for a fixed $n$ large enough and all $\eps$ small enough proves
\eqref{z5.1}.
\end{pf}

\begin{pf*}{Proof of \eqref{z5.4}}
We need a lower bound for
$B(N,t)=\frac{1}2 E_{z_0,\eta_0} [ \mathbf1_{z_t=N}  \eta_t(N-1)
 ]$.
We condition on $\mathcal F_{t-T_1}$ (the canonical filtration) and
denote by
$E_{\bar z, \bar\eta, t-T_1}$ the conditional expectation given
$(\bar z,\bar\eta)$, $\bar z\ne\varnothing$, the configuration at
time $t-T_1$.
The realizations where $z_{t-T_1}=\varnothing$ evidently do not
contribute to $B(N,t)$.

Let $\mathcal D$ denote the event where the rate ${\varepsilon j}/{2}$
clocks at $\pm N$
(attempts to create or remove a particle) never ring in the time
interval $[t-T_1,t]$,
and by $P(\mathcal D)$ its probability. Then
%
\begin{eqnarray}
\label{z5.10} E_{\bar z,\bar\eta, t-T_1} \bigl[ \mathbf1_{z_t=N}
\eta_t(N-1) \bigr] &\ge& E_{\bar z,\bar\eta, t-T_1} \bigl[
\mathbf1_{\mathcal D}  \mathbf 1_{z_t=N} \eta_t(N-1)
\bigr]\nn
\\
&=& P[\mathcal D]  \sum_y q_{T_1}
\bigl(X, (\bar z,y)\bigr) \bar\eta(y)
\\
&=&\mathrm{e}^{-\eps^a j }  \sum_y q_{T_1}
\bigl(X, (\bar z,y)\bigr)\bar\eta(y)\nonumber,
\end{eqnarray}
where $X=(N,N-1)$, $Y=(y_1,y_2)$ and $ q_s(X,Y)$ is the probability
under the stirring process (SSEP)
on $\Lambda_N$ of going from $X$ to $Y$ in a time $s$; the first
equality follows because
the process conditioned on $\mathcal D$ has the law of the stirring
process and the second because
$ P[\mathcal D]=\mathrm{e}^{-\eps j T_1}= \mathrm{e}^{-\eps^a j }$.

Writing $Y=(\bar z,y)$, $Z=(z_1,z_2)$, $Z^0=(z^0_1,z^0_2)$, $z_i\in\La
_N$, $z^0_i\in\La_N$, $i=1,2$:
\[
q_{T_1}(X, Y)=\sum_{Z,Z^0}Q_{T_0}
\bigl(X,X;Z,Z^0\bigr) q_{T_1-T_0}(Z,Y),
\]
where $T_0$ is defined in \eqref{z5.2} and $Q$ refers to the law of
the coupling
between two stirring $(z_1(s),z_2(s))$ and two independent
$(z^0_1(s),z^0_2(s))$ particles as
defined in \cite{DPTVpro} (see Definitions 1 and 4 there in the
particular case of two particles),
with $Q_{T}(\cdot,\cdot)$ denoting the corresponding transition
probabilities in time $T$.
The coupling is such that $z_1(s)=z_1^0(s)$ for all $s\ge0$, and
$z_2(s)$ makes the same jumps as
$z^0_2(s)$ unless $|z_1(s)-z_2(s)|=1$ or one of the involved particles
(independent and stirring)
is at the boundary of $\Lambda_N$ and the other is not.
In particular, if starting at the same pair, independent and stirring
particles move together while
$|z_1(\cdot)-z_2(\cdot)|\ge2$. Moreover, given any $\zeta>0$, the
following estimate is contained in
Theorem 4.5 of \cite{DPTVpro}: for any $n$ there is $c_n$ so that
%
\begin{equation}
\label{z5.11} \sum_{(Z,Z^0)\in\mathcal A^c}Q_{T_0}
\bigl(X,X;Z,Z^0\bigr)\le c_n\eps^n,
\end{equation}
where
%
\begin{equation}
\label{z5.12} \mathcal A=\bigl\{ \bigl(Z,Z^0\bigr):
z_1=z_1^0;\bigl|z_2
-z_2^0\bigr| \le\eps^{-(1-a)/4
-\zeta}\bigr\}.
\end{equation}

\begin{Remark*}
For the case of particles moving in $\mathbb
Z$, this type of estimate has been proven and used
since long ago (see Section~6.6 in \cite{DP}; also Section~3 in \cite{FPSV} or references therein): its rough content is that a pair of
stirring particles can be coupled to a pair of independent random walks
in a way that the first
components coincide, and at time $s$ the second components differ by at
most $s^{1/4 +\delta}$,
except for a set of probability at most $c_k s^{-k}$, as described
above, for any given $\delta>0$.
The restriction to $\Lambda_N$ brings in extra nuisance, as
treated in the proof of Theorem 4.5 of \cite{DPTVpro}.
\end{Remark*}

Let
%
\begin{equation}
\label{z5.13} \mathcal B=\bigl\{Z^0: \bigl|z^0_1
-z^0_2 \bigr| \ge\eps^{-\vfrac{1-a}2 +\zeta}\bigr\}
\end{equation}
so that
%
\begin{equation}
\label{z5.14} q_{T_1}(X, Y)\ge\sum_{(Z,Z^0)\in\mathcal A, Z^0 \in\mathcal
B}Q_{T_0}
\bigl(X,X;Z,Z^0\bigr)q_{T_1-T_0}(Z,Y).
\end{equation}
We write (see \eqref{z5.3})
%
\begin{equation}
\label{z5.15}  \sum_y q_{T_1-T_0}
\bigl(Z, (\bar z,y)\bigr) \bar\eta(y) = p_{T_1-T_0}(z_1, \bar
z)\sum_y p_{T_1-T_0}(z_2, y)
\bar\eta(y) + R(Z),
\end{equation}
where
%
\begin{equation}
\label{z5.15a} R(Z)=\sum_y \bigl[q_{T_1-T_0}
\bigl(Z, (\bar z,y)\bigr)- p_{T_1-T_0}(z_1, \bar
z)p_{T_1-T_0}(z_2, y) \bigr]\bar\eta(y ).
\end{equation}
But if $(Z,Z^0)\in\mathcal A$ and $Z^0 \in\mathcal B$, then $Z\in
\mathcal B':=\{Z\colon|z_1 -z_2 | \ge
\frac{1}2\eps^{-\vfrac{1-a}2 +\zeta}\}$ for all small $\eps$. Also
observe that if we let
\[
\mathcal C= \Bigl\{ \sup_{0\le s\le T_1-T_0} \bigl|z_i(s)-z_i\bigr|
\le (T_1-T_0)^{1/2}\eps^{-\zeta},   i=1,2
\Bigr\}
\]
then whenever $Z\in\mathcal B'$ and $Z(\cdot) \in\mathcal C$,
 we have for $\eps$, $a$, $\zeta$ small enough
\[
\bigl|z_1(s)-z_2(s)\bigr|\ge\tfrac{1}2\eps^{-\vfrac{1-a}2 +\zeta}
- 2 \eps ^{-\vfrac{1-a}4-\zeta}\ge2,\quad\quad 0\le s \le T_1-T_0.
\]
Therefore independent and stirring particles starting from $Z$ can be
coupled to evolve together while in $\mathcal C$, yielding
\[
\mathbb E_Z [\mathbf1_{Z(T_1-T_0)=Y} \mathbf1_{\mathcal C} ]=
\mathbb E^0_Z [\mathbf1_{Z^0(T_1-T_0)=Y}
\mathbf1_{\mathcal
C} ],\quad\quad Z\in\mathcal B',
\]
where $\mathbb E_Z$ and $\mathbb E^0_Z$ ($\mathbb P_Z$ and $\mathbb
P^0_Z$) denote the expectation (law)
relative to the stirring and the independent processes both starting
from $Z$. It follows at once from this and
\eqref{z5.15a} that for $Z\in\mathcal B'$:
\[
\bigl|R(Z)\bigr|\le\mathbb P_Z\bigl[\mathcal C^c\bigr]+\mathbb
P^0_Z\bigl[\mathcal C^c\bigr]\le4 \sup
_{z \in\Lambda_N} \mathcal{P}^0_{z}\Bigl[ \sup
_{0\le s\le T_1-T_0} \bigl|z(s)-z\bigr|> (T_1-T_0)^{1/2}
\eps^{-\zeta}\Bigr],
\]
where at the last inequality we use that under $\mathbb P_Z$ or
$\mathbb P^0_Z$, the components perform simple random walks in $\Lambda
_N$, whose law is written as $\mathcal P^0$.
We then easily see that for each $n$ there exists $c_n$ positive
constant so that
%
\begin{equation}
\label{z5.16} \bigl|R(Z)\bigr| \le c_n\varepsilon^n.
\end{equation}
From \eqref{z5.10}, \eqref{z5.14} and \eqref{z5.16}, we then
get\footnote{Changing the constants $c_n$.}
%
\begin{eqnarray}
\label{z5.17} \nn E_{\bar z, \bar\eta, t-T_1} \bigl[ \mathbf1_{z_t=N}
\eta_t(N-1) \bigr] &\ge& \mathrm{e}^{-\eps^a j } \sum
_{(Z,Z^0)\in\mathcal A, Z^0 \in\mathcal B}Q_{T_0}\bigl(X,X;Z,Z^0\bigr)
\\[-8pt]\\[-8pt]
&&{} \times p_{T_1-T_0}(z_1, \bar z) \sum
_{y\ne\bar z}p_{T_1-T_0}(z_2, y) \bar\eta(y) -
c_n\eps^n.\nonumber
\end{eqnarray}
Letting
%
\begin{equation}
\label{z5.18} \mathcal G=\biggl\{(\bar z, \bar\eta):\bar z\ne\varnothing, \inf
_x \sum_{y\ne\bar z}
p_{T_1-T_0}(x,y)\bar\eta(y) \ge\delta^*\biggr\},
\end{equation}
we can thus write for $\tilde z \neq\varnothing$,
%
\begin{eqnarray}
\label{z5.19} &&E_{\bar z, \bar\eta, t-T_1} \bigl[ \mathbf1_{z_t=N}
\eta_t(N-1) \bigr] \nonumber
\\[-8pt]\\[-8pt]
&&\quad\ge \mathrm{e}^{-\eps^a j } \delta^*
\mathbf1_{\mathcal G}(\bar z,\bar \eta) \sum_{(Z,Z^0)\in\mathcal A, Z^0 \in\mathcal B}Q_{T_0}
\bigl(X,X;Z,Z^0\bigr) p_{T_1-T_0}(z_1, \bar z) -  c_n\eps^n.\nonumber
\end{eqnarray}
But
\begin{eqnarray*}
 && \sum_{(Z,Z^0)\in\mathcal A, Z^0\in\mathcal B}Q_{T_0}
\bigl(X,X;Z,Z^0\bigr) p_{T_1-T_0}(z_1,\bar z)
\\
&&\quad\ge-
Q_{T_0}\bigl(X,X;\mathcal A^c\bigr)
  + \sum_{|z^0_1-z^0_2| \ge\eps^{-\vfrac{1-a}2 +\zeta}}p_{T_0}\bigl(N,
z^0_1\bigr)p_{T_0}\bigl(N-1,
z^0_2\bigr) p_{T_1-T_0}\bigl(z^0_1,
\bar z\bigr),
\end{eqnarray*}
and for any $z^0_1$ and small $\eps$
\[
\sum_{z^0_2: |z^0_1-z^0_2| \ge\eps^{-(1-a)/2
+\zeta}} p_{T_0}\bigl(N-1,z^0_2
\bigr) \ge\frac{1}2,
\]
so that by \eqref{z5.11}
\[
 \sum_{Z,Z^0\in\mathcal A, Z^0 \in\mathcal B}Q_{T_0}
\bigl(X,X;Z,Z^0\bigr) p_{T_1-T_0}(z_1,\bar z) \ge
\frac{1}2 p_{T_1}(N, \bar z) - c_n
\eps^n.
\]
Recalling the definition of $B(N,t)$ and taking the expectation in
\eqref{z5.19} we have
\begin{eqnarray*}
B(N,t)&\ge &\mathrm{e}^{-\eps^a j } \frac{\delta^*}4 \sum
_{z\ne\varnothing} p_{T_1}(N,z)\pi(z,t-T_1)
\\
&&{}- \mathrm{e}^{-\eps^a j } \frac{ \delta^*}{2} P_{z_0,\eta_0}\bigl[ \mathcal
G^c\cap\{z_{t-T_1}\ne\varnothing\}\bigr]- c_n\eps
^n P_{z_0,\eta_0}[z_{t-T_1}\ne\varnothing].
\end{eqnarray*}

In Section~\ref{sec:z6}, we shall prove that
%
\begin{equation}
\label{z5.20} P_{z_0,\eta_0}\bigl[ \mathcal G^c\cap
\{z_{t-T_1}\ne\varnothing\}\bigr] \le c_n \eps^n
P_{z_0,\eta_0}[ z_{t-T_2}\ne\varnothing]
\end{equation}
which will then complete the proof of \eqref{z5.4}.
\end{pf*}

\begin{pf*}{Proof of \eqref{z5.5}}
(The proof given below
uses that the cardinality $K$ of $I_{\pm}$ is 2, for $K>2$ the proof
is similar but more complex.)
By conditioning on $\tilde z_{t-T_1}$, we get
%
\begin{equation}
\label{z5.21} P_{z_0,\eta_0} [ z_t=N ]=\tilde P_{z_0}
[ \tilde z_t=N ] = \tilde E_{z_0} \bigl[
\mathbf1_{\tilde z_{t-T_1}\ne\varnothing
}\tilde P_{t-T_1,\tilde z_{t-T_1}} [\tilde z_t=N ]
\bigr],
\end{equation}
where $\tilde P_{ t-T_1,z'}$ is the law of the auxiliary Markov
process\footnote{Fixed $z_0,\eta_0$ at time $0$ as before.} $\tilde
z_s, s\ge t-T_1$
which starts at time $t-T_1$ from $z'\ne\varnothing$.
Denoting as before by $\mathcal P$ and $\mathcal E$ the law and
expectation of the auxiliary process
with generator $\mathcal{L}^a_t$, that is, when the death part of the
generator is dropped, we have by~\eqref{z3.5},
%
\begin{equation}
\label{z5.22}   \tilde P_{t-T_1,z'} [\tilde z_t=N ] \le
\mathcal P_{t-T_1,z'} [ {\tilde z_{t}= N} ].
\end{equation}
By the integration by parts formula,
\[
\mathcal P_{t-T_1,z'} [ \tilde z_{t}= N ] \le\mathcal
P^0_{T_1}\bigl(z',N\bigr)+ \int
_{t-T_1}^{t} \mathcal P^0_{t-s}(N-1,N)
\frac{\eps j}2 \mathcal P_{t-T_1,z'} [\tilde z_s=N ]\,\mathrm{d}s
+c_k\eps^k,
\]
where $c_k\eps^k$ bounds the contribution of trajectories that visit
$I_-$ and reach $N$ within
time $T_1$ and we used that the rates $a(N,s)$ of extra jumps are
bounded by $\eps j/2$; see \eqref{z3.2}. The random walk probabilities
$\mathcal P^0_{T_1}(z',N)$ and $\mathcal P^0_{t-s}(N-1,N)$
can be computed with the time reverted, yielding
\[
\mathcal P_{t-T_1,z'} [ \tilde z_{t}= N ] \le p_{T_1}
\bigl(N,z'\bigr)+ \int_{t-T_1}^{t}
p_{t-s}(N,N-1) \frac{\eps j}2 \mathcal P_{t-T_1,z'} [\tilde
z_s=N ]\,\mathrm{d}s +c_k\eps^k.
\]
Iterating (and writing $s_0=t$)
%
\begin{eqnarray}
\label{zzzz} && \mathcal P_{t-T_1,z'} [ {z_{T_1}= N} ]\nonumber \\
&&\quad\le\sum
_{n=0}^\infty \biggl( \frac{\eps j}2
\biggr)^n  \int_{t-T_1}^{t} \mathrm{d}s_1
\int_{t-T_1}^{s_1} \mathrm{d}s_2\cdots\int
_{t-T_1}^{s_{n-1}} \mathrm{d}s_n
\\
&&\quad\quad p_{t-s_1}(N,N-1)p_{s_1-s_2}(N, N-1)\cdots \bigl(p_{s_n-(t-T_1)}
\bigl(N, z'\bigr)+c_k \eps^{k} \bigr).\nonumber
\end{eqnarray}
We write the $n$th term of the series as $R_n+R'_n$ where $R_n$ is the
term with $s_n\le t-1$
and $R'_n$ the one with $s_n>t-1$. We start by bounding $R'_n$.
After a change of variables ($s_i\to t-s_i$), calling $\und s=(s_1,\ldots,s_n)$
and $s_0\equiv0$,
%
\begin{eqnarray}\label{z5.23}
R'_n &:=& \biggl( \frac{\eps j}2
\biggr)^n  \int_{[0,T_1]^n, s_n< 1} \Biggl\{\prod
_{i=1}^n \mathbf1_{s_i\ge s_{i-1}}p_{s_i-s_{i-1}}(N,N-1)
\Biggr\} \bigl( p_{T_1-s_n}\bigl(N, z'\bigr)+ c_k
\eps^{k} \bigr)\,\mathrm{d}\und s \nn
\\
&\le& \biggl( \frac{\eps j}2\biggr)^n  \int
_{[0,1]^n} \Biggl\{\prod_{i=1}^n
\mathbf1_{s_i\ge s_{i-1}}\Biggr\} \bigl(p_{T_1-s_n}\bigl(N, z'
\bigr)+c_k \eps^{k} \bigr)\,\mathrm{d}\und s
\\
&\le& \frac{1}{n!} \biggl( \frac
{\eps j}2\biggr)^n \bigl(e
p_{T_1}\bigl(N, z'\bigr)+c_k
\eps^{k} \bigr).  \nn
\end{eqnarray}
To prove the last inequality, we have written
\[
p_{T_1-s_n}\bigl(N, z'\bigr) = \frac{p_{s_n}(N,N)}{p_{s_n}(N,N)}p_{T_1-s_n}
\bigl(N, z'\bigr)\le\frac{p_{T_1}(N, z') }{p_{s_n}(N,N)}
\]
and used $p_{s_n}(N,N)>\mathrm{e}^{-1}$.

To bound $R_n$, we do the same change of variables as above and use the
inequality
\[
p_{s_i-s_{i-1}}(N, N-1)\le\frac{c}{\sqrt{s_i-s_{i-1}}}.
\]
Then
\[
R_n \le\biggl( \frac{\eps j}2\biggr)^n  \int
_{[0,T_1]^n} \mathbf1_{s_n\geq1}  f(\und s) \bigl(
p_{T_1-s_n}\bigl(N, z'\bigr)+c_k
\eps^{k} \bigr)\,\mathrm{d}\und s ,
\]
where
\[
f(\und s) =\mathbf1_{0\equiv s_0\le s_1\le s_2\le\cdots\le s_n\le T_1} \prod_{i=1}^n
\frac{c}{\sqrt{s_i-s_{i-1}}}.
\]
Since $p_{s_n}(N,N)>b /\sqrt{s_n}$ (recall that $s_n\ge1$) we get
\begin{eqnarray*}
  R_n &\le&\biggl( \frac{\eps j}2\biggr)^n\int
_{[0,T_1]^n, s_n\ge1} f(\und s) \biggl(\frac{p_{s_n}(N,N)}{p_{s_n}(N,N)} p_{T_1-s_n}
\bigl(N,z'\bigr)+c_k\eps ^k \biggr)\,\mathrm{d}\und s
\\
&   \le&\biggl( \frac{\eps j}2\biggr)^n \bigl(b^{-1}p_{T_1}
\bigl(N, z'\bigr)+c_k\eps^k \bigr) \int
_{[0,T_1]^n, s_n\ge1} f(\und s)\sqrt{s_n}\,\mathrm{d}\und s.
\end{eqnarray*}
We change variables: $s_i\to T_1 s_i$ and get, using Lemma 5.2 of
\cite{DPTVjsp},
\begin{eqnarray*}
\int_{[0,T_1]^n, s_n\ge1} f(\und s)\sqrt{s_n}\,\mathrm{d}\und s &\le&
T_1^{(n+1)/2} \int_{[0,1]^n} f(\und s)
\sqrt{s_n}\,\mathrm{d}\und s
\\
& \le& T_1^{(n+1)/2} \int_{[0,1]^n} f(\und s)\,\mathrm{d}\und s
\\
&\le& C^n \mathrm{e}^{-\ifrac{n}2[\log\ifrac{n}2 -1]} \eps^{-\ifrac{1}2(n+1)
+\ifrac{a} 2(n+1)}.
\end{eqnarray*}
Thus,
%
\begin{eqnarray}
R_n &\le& \biggl(\frac{C j}{2}\biggr)^n
\mathrm{e}^{-\ifrac{n}2[\log\ifrac{n}2 -1]} \eps ^{\ifrac{1}2(n-1)+\ifrac{a} 2(n+1)} \bigl(b^{-1}p_{T_1}
\bigl(N, z'\bigr) +c_k\eps^k \bigr).
\label{z5.24}
\end{eqnarray}

Putting together the estimates \eqref{z5.23} and \eqref{z5.24}, we
can bound the sum of $R_n +R'_n$ over $n$ in \eqref{zzzz} (by
convergent series). It follows that positive constants $\tilde c$ and
$c_k$ can be found so that for all $z' \in\Lambda_N$ and all $k$
\[
\mathcal P_{t-T_1,z'} [ \tilde z_{t}= N ] \le\tilde c
p_{T_1}\bigl(N,z'\bigr) + c_k
\eps^k
\]
for all $\eps$ small. Now combining this into \eqref{z5.21}, and since
$\tilde P_{z_0}[\tilde z_{t-T_1} \ne\varnothing] \le\tilde
P_{z_0}[\tilde z_{t-T_2} \ne\varnothing]$ we
have \eqref{z5.5}.
\end{pf*}

\begin{pf*}{Proof of \eqref{z5.6}}
Let $t\ge T_2 := \kappa\eps^{-2}$,
then analogously to \eqref{z3.5},
%
\begin{equation}
\pi(N,t)\equiv\tilde P_{z_0}[\tilde z_t =N] = \tilde{ E
}_{z_0} \bigl[\mathbf1_{\tilde z_{ t-T_2} \ne\varnothing} {\mathcal E}_{t-T_2,
\tilde z_{ t-T_2}}
\bigl[\mathrm{e}^{ - \int_{t-T_2}^td( z_s,s)\,\mathrm{d}s }\mathbf 1_{ z_{ t}=N} \bigr] \bigr] \label{z5.25}
\end{equation}
with ${\mathcal E}_{t,x}$ as defined in the beginning of Section~\ref{sec:4}.

We
denote by $\mathcal E'_{N}$ the expectation with
respect to the time-backward process, $z'_s$, $s\in[0,T_2]$, which
starts at time $0$ from $N$ and is a simple random walk
with additional jump intensity $a(\pm N,t-s)$ for the jump $\pm
(N-1)\to\pm N$ at time $s$. We then have
%
\begin{eqnarray}
\label{z5.26} \pi(N,t) &=& \mathcal E'_{N} \biggl[\pi
\bigl(z'_{T_2}, t-T_2\bigr) \exp\biggl\{ -
\int_0^{T_2}  d \bigl(z'_s,t-s
\bigr)\,\mathrm{d}s \biggr\} \biggr]\nonumber
\\
&\ge&\mathrm{e}^{-\eps j} \mathcal E'_{N} \biggl[\pi
\bigl(z'_{T_2}, t-T_2\bigr) \exp\biggl\{ -
\int_1^{T_2-1}  d \bigl(z'_s,t-s
\bigr)\,\mathrm{d}s \biggr\} \biggr] \nn
\\
&\ge& \mathrm{e}^{-\eps j} \mathcal E'_{N} \biggl[\pi
\bigl(z'_{T_2}, t-T_2\bigr)
\mathbf1_{z'_{1}=N-2} \exp\biggl\{ - \int_1^{T_2-1}  d \bigl(z'_s,t-s\bigr)\,\mathrm{d}s \biggr\} \biggr]
\\
&\ge& \mathrm{e}^{-\eps j} \alpha\sum_{|x|\le N-2}
\pi(x,t-T_2)     \alpha'{\mathcal P}^0_{N-2}
\Bigl[x_{T_2-2}=x,\sup_{s\in[0,T_2-2]}|x_s|<N-1
\Bigr] \nn
\\
&&{}+ \mathrm{e}^{-\eps j} \alpha
  \sum_{\sigma=\pm,x\in I_{\sigma}}\pi(x,t-T_2)
    \alpha'' \mathcal{ P}^0_{N-2}
  \Bigl[x_{T_2-2}=\si (N-2), \sup_{ s\in [0,T_2-2]}|x_s|<N-1 \Bigr],\nn
\end{eqnarray}
where ${\mathcal P}^0_{x}$ is the law of the random walk $x_s$ with no
extra jumps (just
a simple random walk on $\Lambda_N$ starting at $x$) and
\begin{eqnarray*}
\alpha&=& \mathcal P'_N\bigl[z'_1=N-2
\bigr] >0, \quad\quad\alpha'= \min_{|x|\le N-2} {\mathcal
P}^0_N[x_{T_2-1}= x|x_{T_2-2}=x] >0,
\\
\alpha''&=&\min_{\sigma\in\{-1,1\}}\min_{x\in\{N-1,N\}}
\mathcal{ P}'_N\bigl[z'_{T_2}= \sigma x|z'_{T_2-1}=\si(N-2)\bigr] >0.
\end{eqnarray*}
We thus need to bound from below
the probability of the event
$\{x_{T_2-2}=x, |x_s|\le N-2, s\in[0, T_2-2]\}$
uniformly in $|x| \le N-2$. The basic idea is to reduce to a single
time estimate; indeed, the
condition $|x_s|\le N-2, s\in[0, T_2-2]$ can be dropped provided we
study the
process on the whole $\mathbb Z$ and take as initial condition the
antisymmetric datum, which
is obtained by assigning a weight $\pm1$ to the images of $x$ under
reflections around
$\pm(N-1)$. The details are given in the \hyperref[app]{Appendix}. To have control of
the plus and minus contributions, it is convenient to reduce
to small time intervals; moreover, the analysis will
distinguish the case where $x$ is ``close'' to $\pm N$ and
when it is not. Closeness here means that
$N-|x| \le N/100$ (the choice $1/100$
is just for the sake of concreteness, any ``small'' number would work
as well).

Let us now be more specific.
We split $T_2 -2= m\tau\eps^{-2}$, $m$ an integer and $\tau>0$ small
enough, and write
\begin{eqnarray*}
&& \mathcal{P}^0_{N-2} \Bigl[x_{m\eps^{-2}\tau}=x;\sup
_{ s\in[0,T_2-2]}|x_s|<N-1 \Bigr]
\\
&&\quad \ge\mathcal{P}^0_{N-2} \Biggl[ \bigcap
_{i=1}^{m-1} \Bigl\{ \sup_{ s\in[i-1,i]\eps^{-2}\tau
}|x_s|<N-1;
|x_{i\eps^{-2}\tau}| \le N/100\Bigr\}
\\
&&\quad\hphantom{\ge\mathcal{P}^0_{N-2} \Biggl[ \bigcap
_{i=1}^{m-1}}{} \cap\Bigl\{\sup_{s\in[m-1,m]\eps^{-2}\tau}|x_s|<N-1;
x_{m\eps
^{-2}\tau}=x\Bigr\} \Biggr].
\end{eqnarray*}
In the \hyperref[app]{Appendix}, we shall prove that for $\tau$ small enough there is
$c$ so that for all $\eps$ ($N=\eps^{-1}$),   the following bounds hold:
%
\begin{eqnarray}
\label{5.27} \mathcal{P}^0_{N-2} \Bigl[|x_{\eps^{-2}\tau}|
\le N/100; \sup_{ s\in[0,\eps^{-2}\tau
]} |x_s|<N-1 \Bigr] &\ge& c \eps,
\\
\label{5.28} \inf_{|x|\le N/100} \mathcal{P}^0_{x}
\Bigl[|x_{\eps^{-2}\tau}| \le N/100; \sup_{ s\in[0,\eps^{-2}\tau
]}
|x_s|<N-1 \Bigr] &\ge& c,
\\
\label{5.29} \inf_{|x|\le N/100} \inf_{|x'| \le N 99/100}
\mathcal{P}^0_{x} \Bigl[x_{\eps^{-2}\tau}=x';
\sup_{ s\in[0,\eps^{-2}\tau]} |x_s|<N-1 \Bigr] &\ge& c \eps,
\\
\label{5.30} \inf_{|x|\le N/100} \inf_{ N 99/100 \le|x'| \le N-2}
\mathcal{P}^0_{x} \Bigl[x_{\eps^{-2}\tau}
=x'; \sup_{ s\in[0,\eps^{-2}\tau]} |x_s|<N-1 \Bigr]
&\ge& c \eps^2.
\end{eqnarray}
The above bounds together with \eqref{z5.26}
prove \eqref{z5.6}.
\end{pf*}

\section{Proof of \texorpdfstring{\protect\eqref{z5.20}}{(5.21)}}
\label{sec:z6}
For any $(z,\eta)$, we define the configurations $\eta^{(1)}$ and
$\eta^{(2)}$ in $\{0,1\}^{\La_N}$ as follows:
If $z\ne\varnothing$, then $\eta^{(1)}(x)=\eta^{(2)}(x)=\eta(x)$ for
any $x\in\La_N\setminus{z}$,
and $\eta^{(1)}(z)=1$, $\eta^{(2)}(z)=0$. If $z=\varnothing$, then
$\eta^{(1)}=\eta^{(2)}=\eta$.

If $(z_t,\eta_t)_{t\ge0}$ is the process defined in Section~\ref{sec:z2}, we can see that $(\eta^{(2)}_t)_{t\ge0}$ has the
law of the process introduced in \cite{DPTVjsp} that we are here
calling DPTV for simplicity (as well as $(\eta^{(1)}_t)_{t\ge0}$,
though such a property will not be used in the following). Details
can be found in \cite{DPTVgap}.

For any $x\in\La_N$, we introduce the function
$ A_x(\eta)$, $\eta\in\{0,1\}^{\La_N}$, by setting
%
\begin{equation}
\label{6.2} A_x(\eta) := \sum_{y}
p_{T_1-T_0}(x,y)\eta(y), \quad\quad\eta\in\{0,1\}^{\La_N}.
\end{equation}
Then, recalling that $\mathcal G$ has been defined in \eqref{z5.18}
and writing $\tau:= t- T_1$,
the left-hand side of \eqref{z5.20} is equal to
\[
P_{z_0,\eta_0} \Bigl[ z_\tau\ne\varnothing, \inf
_x A_x\bigl(\eta ^{(2)}_\tau
\bigr) \le\delta^* \Bigr] \le E_{z_0,\eta_0} \Bigl[ \mathbf1_{z_{t-T_2} \ne\varnothing}
P_{z_{t-T_2},\eta_{t-T_2}}\Bigl[ \inf_x A_x\bigl(
\eta^{(2)}_\tau\bigr) \le\delta^* \Bigr] \Bigr]
\]
which is bounded from above by
\[
\tilde P_{z_0} [\tilde z_{t-T_2} \ne\varnothing] \sup
_{\eta\in\{0,1\}
^{\La_N}} \mathbf P_{\eta}\Bigl[   \inf
_x A_x(\eta_{T_2-T_1}) < \delta^*
\Bigr],
\]
where $\mathbf P_{\eta}$ is the law of the DPTV process starting from
$\eta$ at time $0$. We thus need to prove that
\[
\sup_{\eta\in\{0,1\}^{\La_N}} \mathbf P_{\eta}\Bigl[  \inf
_x A_x(\eta_{T_2-T_1}) < \delta^* \Bigr]
\le c_n \eps^n.
\]
Since the evolution preserves the coordinate-wise order in $\{0,1\}
^{\Lambda_N}$ (see \cite{DPTVpro}) and $\inf_x A_x(\eta) $ is a
non-decreasing
function of $\eta$, it suffices to show that\vspace*{-1pt}
%
\begin{equation}
\label{6.3} \mathbf P_{\mathbf0}\Bigl[  \inf_x
A_x(\eta_{T_2-T_1}) < \delta^* \Bigr] \le c_n
\eps^n,
\end{equation}
with $\mathbf{0}$ the configuration with $\eta(x)=0$ for all $x$.

In \cite{DPTVpro}, it is proved that there is $\tau^*>0$ (independent
of $N$) so that if $t \in[N^2,\tau^*N^2\log N]$
then for any $n$ there is $c_n$ so that\vspace*{-1pt}
%
\begin{equation}
\label{00app.5} \mathbf P_{\mathbf0} \Bigl[\inf_x
\bigl|A_x(\eta_t) - A_x\bigl(\ga(\cdot ,t)
\bigr)\bigr| \ge\eps^{1/4} \Bigr] \le c_n \eps^n,
\end{equation}
where $\ga(y,t)= \rho(\eps y, \eps^2 t)$ and $\rho(r,t)$, $r\in
[-1,1]$, $t\ge0$, is the solution of the hydrodynamic equation
for the DPTV system starting from $\rho(r,0)\equiv0$. In \cite{DPTVjsp2},
it is proved that\vspace*{-1pt}
%
\begin{equation}
\label{00app.6} \lim_{t\to\infty} \sup_{|r|\le1} \bigl|
\rho(r,t)-\rho^{\rm st}(r)\bigr| = 0
\end{equation}
and that $\rho^{\rm st}(r)$ is an increasing function (linear with
positive slope) with
$\rho^{\rm st}(-1)>0$. Thus, there is $\kappa>0$ independent of $N$
so that for all $N$ large enough\vspace*{-1pt}
%
\begin{equation}
\label{00app.7} \mathbf P_0 \biggl[\inf_x
A_x(\eta_s)\ge\frac{\rho^{\rm st}(-1)}2 \biggr]
\ge1-c_n \eps^n, \quad\quad\frac{\kappa}{2} N^2\le s
\le\kappa N^2
\end{equation}
which implies \eqref{6.3}, provided $\delta^* < \rho^{\rm st}(-1)/2$
and $T_2 = \kappa N^2$.\vspace*{-2pt}

\begin{appendix}
\section*{Appendix}\label{app}\setcounter{equation}{0}

We now prove the bounds \eqref{5.27}--\eqref{5.30}. The key point is
the identity below for
the transition probabilities for the simple random walk in an interval,
absorbed at the
boundaries. (The proof follows the same argument as that given for
the Brownian motion case; see, e.g., Proposition 8.10 in Chapter~2 of
\cite{KS}.)
Let $L=N-1\ge2$,\vspace*{-1pt}
%
\begin{eqnarray}
\label{A.1} &&\mathcal{P}^0_{x} \bigl[{x_{t}
=y}; {|x_s|<L, \forall s\in\bigl[0,\eps^{-2}\tau\bigr]}
\bigr]\nonumber\\[-8.5pt]\\[-8.5pt]
&&\quad = \sum_{k\in
\mathbb Z} \bigl[p_t(4kL+y-x)
-p_t(4kL-2L-y-x)\bigr],\nonumber
\end{eqnarray}
where $x$ and $y$ in \eqref{A.1} are in $[-L+1,L-1]$, and $p_t(z)$ is
the probability for a
simple random walk on $\mathbb Z$ starting from $0$ to be at $z$ at
time $t$, for $z \in\mathbb{Z}$.

Writing $z=L-y$ and $w=L-x$, rearranging the sum, and using the
symmetry of $p_t(\cdot)$, we rewrite \eqref{A.1} as\vspace*{-1pt}
%
\begin{eqnarray}
\label{A.2} && \mathcal{P}^0_{x}
\bigl[{x_{t} =y}; {|x_s|<L, \forall s\in\bigl[0,
\eps^{-2}\tau\bigr]} \bigr] \nonumber
\\[-0.5pt]
&&\quad= p_t (z-w ) -
p_t(z+w)+\sum_{k=1}^\infty \bigl(
\bigl[p_t (4kL-z+w )-p_t (4kL-z-w )\bigr]
\\[-0.5pt]
&&\quad\quad\hphantom{p_t (z-w ) -
p_t(z+w)+\sum_{k=1}^\infty \bigl(}{}  -
\bigl[p_t (4kL+z+w )-p_t (4kL+z-w )\bigr] \bigr).\nn
\end{eqnarray}

To prove \eqref{5.27} (where $x=N-2$), we take $w=1$ in \eqref{A.2}
and get (recall $L=N-1$):
%
\begin{eqnarray}
\label{A.3}  &&\mathcal{P}^0_{N-2}
\bigl[{x_{t} =L-z}; {|x_s|<N-1, \forall s\in[0,t]} \bigr]\nonumber\\
&&\quad\ge p_t(z-1)-p_t(z+1)\nn
\\[-8pt]\\[-8pt]
&&\quad\quad{} - \sum_{1\le k\le\eps^{-b}}\sum_{\si=\pm1}
\bigl\llvert p_t(4kL +\si z-1)-p_t(4kL+\si z+1)\bigr
\rrvert \nonumber
\\
&&\quad\quad{} - 2\sum_{|u|\ge N\eps^{-b}/2} p_t(u),\nn
\end{eqnarray}
$b>0$ a small constant.

Given $b$ and $\tau$ positive constants (independent of $\eps$),
since $t=\tau\eps^{-2}$ there is $c>0$ so that for all $\eps$ small
enough (and, say, all $\tau\in(0,1]$)
%
\begin{equation}
\label{A.4} \sum_{|u|\ge N\eps^{-b}/2} p_{t}(u) \le
\mathrm{e}^{- c\eps^{-2b}},
\end{equation}
by simple tail estimate for the random walk on $\mathbb Z$.

We prove that $p_t(z-1)-p_t(z+1)$ is bounded from below
proportionally to $\eps^2$, so that the last sum in \eqref{A.3} will
be negligible with
respect to the first.
The other terms on the right-hand side
of \eqref{A.3} are bounded in the following proposition, and using the
smallness of
$\tau$ we see that their sum over $1\le k\le\eps^{-b}$
is a small fraction of the first term on the right-hand side of \eqref{A.3},
from which \eqref{5.27} will follow.

\begin{prop}
Recalling that $N\equiv\eps^{-1}$, $t\equiv\eps^{-2}\tau$,
there are positive constants $c$, $C$ and $b$ such that for every $\tau
$, the following holds
for all $\eps$ small enough:
\begin{itemize}
\item  When $N/2< y < 2 N$,
%
\begin{equation}
\label{1} p_t(y)-p_t(y+2)\geq\frac{\eps^{2}}{\sqrt{2\uppi\tau}}
\mathrm{e}^{-(\eps y)^2/2\tau} \frac{1}{4\tau}(1-c\eps).
\end{equation}

\item
When $N/2 <y < N\eps^{-b}$,
%
\begin{equation}
\label{2} p_t(y)-p_t(y+2)\leq\frac{\eps^{2}}{\sqrt{2\uppi\tau}}
\mathrm{e}^{-(\eps y)^2/2\tau} \frac{8 \eps y}{\tau}(1+c\eps).
\end{equation}
\end{itemize}
\end{prop}

\begin{pf}
We
have
\[
p_t(y) = \mathrm{e}^{-t}{\sum_n}^*
\biggl(\frac{1}2\biggr)^n\frac{t^n}{n!} \pmatrix{ n
\cr
m},\quad\quad
y= 2m-n,
\]
where $\sum_n^*$ means that $n$ runs over either the odd or the even
integers of $\mathbb Z$ according to whether
$y$ is odd or, respectively, even.
$n$ is the total number of jumps, $m$ the number of jumps to the right
so that $m-(n-m)=y$.

We start by proving \eqref{1}.
For every pair $y$ and $y':=y+2$, let $m$ and $m'$ be the number of the
corresponding jumps
to the right, so that $m'=m+1$.
Then
%
\begin{equation}
\label{5} \pmatrix{ n
\cr
m}-\pmatrix{ n
\cr
m'}=\pmatrix{ n
\cr
m}\biggl(1-\frac{n-m}{m+1}\biggr)= \pmatrix{ n
\cr
m}\frac{y+1}{m+1}.
\end{equation}
We bound $m = (n+y)/2 \le t$, which is valid when $n \le2t - 2N$. Thus,
\begin{eqnarray*}
p_t(y)-p_t(y+2) & \ge& \frac{N}{2( t+1)}
\mathrm{e}^{-t}\sum_{n \le2t -
2N}\biggl(\frac{1}2
\biggr)^n \frac{t^n}{n!} \pmatrix{ n
\cr
m}
\\
& \ge& \frac{N}{4 t} \mathrm{e}^{-t} \biggl[ \sum
_{n\geq1}\biggl(\frac{1}2\biggr)^n
\frac
{t^n}{n!} \pmatrix{ n
\cr
m} -\sum_{n > 2t - 2N}
\biggl(\frac{1}2\biggr)^n \frac{t^n}{n!} \pmatrix{ n
\cr
m}
\biggr]
\end{eqnarray*}
and \eqref{1} then follows from the local limit theorem (\cite{lawler}, page~58,
Theorem 2.5.6), after observing that the sum over
$n > 2t - 2N$ is exponentially small in $t$.

To prove \eqref{2}, we proceed similarly. Since we want an upper
bound, we
write $m+1 \ge n/2$, getting\vspace*{-0.5pt}
\[
p_t(y)-p_t(y+2) \le\frac{y+1}{t/4}\mathrm{e}^{-t}
\sum_{n \ge t/2}\biggl(\frac{1}2
\biggr)^n \frac{t^n}{n!} \pmatrix{ n
\cr
m}.
\]
As before, \eqref{2} is again a consequence of the local limit
theorem, and the large deviation estimate on
the number of jumps for the set $n<t/2$.\vspace*{-1pt}
\end{pf}

\begin{pf*}{Proof of \eqref{5.27}}
By \eqref{A.3} and \eqref{A.4},
using the above proposition,
\begin{eqnarray*}
  &&\mathcal{P}^0_{N-2} \bigl[{x_{t} =y};
{|x_s|<L,\mbox{ for all } s\in[0,t]} \bigr]
\\
&&\quad\ge\frac{\eps^{2}}{\sqrt{2\uppi\tau}}
\mathrm{e}^{-(\eps z)^2/2\tau}
\frac{1}{4\tau}(1-c\eps)
\\
&&\quad\quad{} - 2\sum_{1\le k\le\eps^{-b}} \frac{\eps^{2}}{\sqrt{2\uppi\tau}}
\mathrm{e}^{-( [4k(1-\eps)-\eps z -\eps])^2/2\tau} \frac{8\eps(4k+2)}{\tau
}(1+c\eps)
\\
&&\quad\quad {}- 2 \mathrm{e}^{-c\eps^{-2b}},\quad\quad \mbox{ where } z =N-1-y.
\end{eqnarray*}
If $\tau>0$ is sufficiently small, then for all $\eps$ small enough\vspace*{-0.5pt}
\[
\mathcal{ P}^0_{N-2} \bigl[{x_{t} =y};
{|x_s|<L,\mbox{ for all } s\in[0,t]} \bigr] \ge\frac{\eps^{2}}{\sqrt{2\uppi\tau}}
 \mathrm{e}^{-(\eps z)^2/2\tau}
\frac{1}{8\tau}
\]
and \eqref{5.27} is proved.
\end{pf*}

To prove \eqref{5.28} and \eqref{5.29}, we use again \eqref{A.1} and bound
%
\begin{equation}
\label{A.8} \mathcal{P}^0_{x} \Bigl[{x_{t}
=y}; \sup_{s\in[0,\eps^{-2}\tau]}|x_s|<L \Bigr] \ge
p_t(y-x) - \sum_{n \in\mathbb{Z}, n\neq0} p_t(
y_n-x)- \sum_{n \in
\mathbb{Z}}p_t
\bigl(y'_n-x\bigr),
\end{equation}
where $y_n=y +4nL$ and $y'_n=-y-2L +4nL$, which we may rewrite as
%
\begin{equation}
\label{A.8'} \mathcal{P}^0_{x} \Bigl[{x_{t}
=y}; \sup_{s\in[0,\eps^{-2}\tau]} |x_s|<L \Bigr] \ge
p_t(\tilde y_0-x) - \sum_{n \in\mathbb{Z}, n\neq0}
p_t( \tilde y_n-x),
\end{equation}
where $\tilde y_0=y$ and the points $\tilde y_n$ stay at distance at
least $aN$ from each other.
As before, we then may bound from below the right-hand side by
\[
p_t(y-x) - \sum_{1\le|n|\le N\eps^{-b}} p_t(
\tilde y_n-x) - \sum_{|z| \ge N\eps^{-b}}
p_t(z)
\]
and \eqref{5.28} and \eqref{5.29} follow using the local limit
theorem and large deviations as before.

\begin{pf*}{Proof of \eqref{5.30}}
We use the equality
%
\begin{equation}
\label{A.9} \mathcal{P}^0_{x} \Bigl[{x_{\eps^{-2}\tau}
=x'}; \sup_{s\in[0,\eps^{-2}\tau]} |x_s|<L \Bigr]=
\mathcal{P}^0_{x'} \Bigl[x_{\eps^{-2}\tau} =x;\sup
_{s\in[0,\eps^{-2}\tau]}|x_s|<N-1 \Bigr]
\end{equation}
recalling that $|x|\le N/100$ and $ N 99/100 \le|x'| \le N-2$; we thus
need to bound
from below the right-hand side of \eqref{A.9} by $c\eps^2$ with $c>0$
independent of
$x$ and $x'$ when they vary in the above
sets.

We thus use \eqref{A.2} with $x\to x'$ and $y\to x$, so that on the
right-hand side
we must read $z=L - x$ and $w=L-x'$. Observe that $N-1-\frac{N}{100}
\le z\le N-1+ N/100$ and $w \in[1,\frac{N}{100}-1] \cup[2N-1-\frac
{N}{100},2N-3]$.
To have the same structure as in \eqref{A.3}, we write
\[
p_t (z-w ) - p_t(z+w) = \sum
_{i=1}^{w} \bigl[p_t(z-y_i-1)
-p_t(z-y_i+1)\bigr],
\]
where $y_i= w + (2i-1), 1\le i\le w$
with the analogous decomposition for $p_t (z'-w ) - p_t(z'+w)$ with
$z'=4kL \pm z$. Then
%
\begin{eqnarray}
\label{A.222}  && \mathcal{P}^0_{x'}
\bigl[{x_{t} =x}; {|x_s|<L,\mbox{ for all } s\in\bigl[0,
\eps^{-2}\tau\bigr]} \bigr] \nonumber
\\
&&\quad\ge\sum_{i=1}^w
\biggl( \bigl[p_t (z-y_i-1) - p_t(z-y_i+1)
\bigr]\nn
\\[-8pt]\\[-8pt]
&&\quad\quad\hphantom{ \sum_{i=1}^w
\biggl(}{} -\sum_{1\le k \le N\eps^{-b}}\sum_{\si=\pm1}
\bigl|p_t \bigl(4kL-\si(z-y)-1 \bigr)-p_t \bigl(4kL-\si(z-
y)+1 \bigr)\bigr|\quad\quad\nonumber
\\
&&\quad\quad\hphantom{ \sum_{i=1}^w
\biggl(}{} - 2\sum_{|u|\ge N\eps^{-b}/2} p_t(u) \biggr)\nn
\end{eqnarray}
and for each $y_i$ we have the same bound as before, hence \eqref{5.30}.
\end{pf*}

\end{appendix}

\section*{Acknowledgements}
The research has been partially supported by PRIN 2009 (2009TA2595-002).
The research of Dimitrios Tsagkarogiannis is partially supported by the FP7-REGPOT-2009-1 project
``Archimedes Center for Modeling, Analysis and Computation'' (under
grant agreement no. 245749).
Maria Eulalia Vares is partially supported by CNPq grants PQ 304217/2011-5 and 474233/2012-0.

The authors are grateful to the referees for their careful reading.
Their comments helped to improve the paper.


%

\printhistory
\end{document}
